\documentclass[reqno]{amsart}
\usepackage{amsfonts}
\usepackage{latexsym}
\author{Hur\c{s}{\.i}t \"{O}ns{\.i}per}
\address{Department of Mathematics\\
Middle East Technical University\\
06531 Ankara, Turkey}
\email{hursit@@rorqual.cc.metu.edu.tr}
\author{S{\.i}nan Sert\"{o}z}
\address{Department of Mathematics\\
Bilkent University\\
06533 Ankara, Turkey}
\email{sertoz@@fen.bilkent.edu.tr}
\subjclass{Primary: 14E15; Secondary: 14J10, 14J15, 14J29}
\thanks{To appear in {\em Archiv der Mathematik.}}
\title[On Degenerations of Fibre Spaces of Curves of Genus $\geq
2$]{On Degenerations of Fibre Spaces \\ of Curves of Genus $\geq 2$}
\theoremstyle{plain}
\newtheorem{theorem}{Theorem}
\newtheorem{lemma}{Lemma}
\newtheorem{corollary}{Corollary}
\begin{document}
\maketitle
\begin{abstract} 
In this note, we show that for surfaces admitting
suitable fibrations, any given degeneration ${\mathcal X}/\Delta$ is
bimeromorphic to a fiber space over $\Delta$ and we apply this
result to the study of the degenerate fiber.\\
\end{abstract}

This note is concerned with the problem of studying the
degenerations of fibered surfaces via the degenerations of the base
curve and the fibers. We consider a surface of general type $ X $
admitting a fibration $X \rightarrow S$ with base genus $\ge 2$ and
show that for a weakly projective degeneration $\Pi : {\mathcal X}
\rightarrow \Delta$ of such a surface, satisfying a mild condition
on monodromy, the components of the singular fiber are of the
following types:\\ (i) a rational surface, (ii) a ruled surface, or
(iii) a surface fibered over a curve obtained from a degeneration of
the base curve $S$.\\

Throughout the paper we work over $\mathbb{C}$, and adapt the
following notation :\\ $\Delta = \{z \in \mathbb{C} : \vert z \vert
< 1\}$.\\ $\Delta^{*} = \Delta - {0}$.\\ X is a compact surface
admitting a fibration $ X \rightarrow S$ with general fiber $F$ and
genera $g(S), g(F) \ge 2$.\\ $\chi(?)$ denotes the holomorphic Euler
characteristic of ?.\\ $K_{?}$ is the canonical class of ?.\\
$p_{g}(?), q(?)$ denote the geometric genus and the irregularity of
?, respectively.\\ ${\mathcal X}$ is a smooth threefold, $\Pi :
{\mathcal X}
\rightarrow \Delta$ is a flat proper holomorphic map with ${\mathcal
X}_{t_{0}} = X$ for some $t_{0} \in \Delta^{*}$ and the only
singular fiber ${\mathcal X}_{0}$ is a divisor with normal crossings
and smooth components.\\

We further assume that (i) each fiber ${\mathcal X}_{t}, t \in
\Delta^{*}$ admits a fibration ${\mathcal X}_{t} \rightarrow S_{t}$
of the same type as the fibration on $X$ and (ii) the monodromy
action on $H^{1}({\mathcal X}_{t}, \mathbb{C})$ leaves the image of
$H^{1}(S_{t},
\mathbb{C})$ invariant (cf. Lemma 3 and the discussion preceding Lemma
3, for some examples of degenerations satisfying these hypothesis).
We recall that under these assumptions we have the following result
([4], Proposition 1 and Proposition 3).\\ 
						  
\begin{theorem} Possibly after restricting to a smaller disk
$\Delta'$ around $0 \in \Delta$, we can find a degeneration $\Pi' :
{\mathcal X}' \rightarrow \Delta'$ and a relative curve $p : {\mathcal Y}
\rightarrow \Delta'$ such that :\\ (i) we have a bimeromorphic map
${\mathcal X}'/\Delta' \rightarrow {\mathcal X}/\Delta'$ which is an
isomorphism over ${\Delta'}^{*}$,\\ (ii) $\Pi'$ factors through $p$,
${\mathcal X}^{'} \rightarrow {\mathcal Y}$ is surjective and\\
(iii) for each $t \in \Delta^{'*}$, ${\mathcal X}_{t}^{'}
\rightarrow {\mathcal Y}_{t}$ is the fibration ${\mathcal X}_{t}
\rightarrow S_{t}$.
\end{theorem}

Theorem 1 reduces the problem of understanding the
structure of the given degeneration to the study of the
degenerations of the curves $S_{t}$ and of the fibres of $X
\rightarrow S$. In this direction, without any projectivity
assumptions on ${\mathcal X}/\Delta$, we have\\

\begin{lemma} If a component $X_{i}$ of ${\mathcal X}_{0}^{'}$
maps onto a component $Y_{i}$ of ${\mathcal Y}_{0}$, then it is
algebraic.
\end{lemma}

\noindent {\it Proof :} 
If $X_{i} \rightarrow Y_{i}$ is a fibration of fiber genus 0 or $\ge
2$ then clearly $X_{i}$ is algebraic. So we assume that $X_{i}
\rightarrow Y_{i}$ is an elliptic fibration. Since a surface with
algebraic dimension = 0 has only a finite number of curves, $X_{i}$
has algebraic dimension = 1 or 2.  As the fibers of ${\mathcal
X}^{'}_{0} \rightarrow {\mathcal Y}_{0}$ are curves of genus $\ge 2$
(except at a finite number of points where the map is not flat), we
see that on $X_{i}$ lie other curves obtained from intersection with
the components of ${\mathcal X}^{'}_{0}$ containing the rest of the
fibers of ${\mathcal X}_{0}^{'} \rightarrow {\mathcal Y}_{0}$ over
$Y_{i}$.  Clearly, such curves will not lie in the fibers of the
elliptic fibration on $X_{i}$ and this is impossible unless $X_{i}$
has algebraic dimension 2. This proves the lemma. \hfill $\Box$\\

\begin{lemma} Let $X_{j}$ be a component of ${\mathcal X}^{'}_{0}$
mapped to a point $p \in {\mathcal Y}_{0}$.  If $X_{j}$ intersects some
component $X_{i}$ as in Lemma 1 along a smooth curve $C$ with no
triple points, then $X_{j}$ is algebraic.
\end{lemma}

\noindent {\it Proof :}
Let $C_{i} (resp. C_{j})$ denote the curve $C$ on $X_{i} (resp.
X_{j})$.\\ If $C_{i}$ is the fiber of $X_{i}$ over $p$, then
$C_{i}^{2} = 0$. Therefore, as there are no triple points on $C$ we
have $C_{j}^{2} = - C_{i}^{2} = 0$. Moreover $g(C) = g(F) \ge 2$ and
using the adjunction formula on $X_{j}$, we get $K_{j}.C_{j} +
C^{2}_{j} = 2(genus(C) - 1) \ge 2$. Hence $(K_{j} + nC_{j})^{2} > 0$
for large enough n and therefore $X_{j}$ is algebraic.\\ On the
other hand, if $C_{i}$ is a component of the fiber over $p$, then
$C_{i}^{2} < 0$ and the equality $C_{j}^{2} = - C_{i}^{2}$ gives
$C_{j}^{2} > 0$, again proving the algebraicity of $X_{j}$.
\hfill $\Box$\\

These two lemmata clearly fall short of proving the algebraicity of
all components of the singular fiber (cf. [5], conjecture on p. 83).
However, combining Lemma 1 with the flatification technique of ([2])
we get \\

\begin{theorem} ${\mathcal X}/\Delta^{'}$ is bimeromorphic to a
degeneration ${\mathcal X}^{''}/\Delta^{'}$ in which all components
of the singular fiber ${\mathcal X}_{0}^{''}$ are algebraic.
\end{theorem}

\noindent {\it Proof :} To prove this result we first remove those components
of ${\mathcal X}^{'}_{0}$ where the map ${\mathcal X}^{'} \rightarrow
{\mathcal Y}$ fails to be flat. For this purpose, we will apply
flatification as described in ([2]). More precisely, we blow up
points $p_{1},..,p_{k} \in {\mathcal Y}_{0}$ over which our map is
not flat, to get a new relative curve ${\mathcal Y}^{'}/\Delta^{'}$
and then in the complex space ${\mathcal X}^{*} = {\mathcal X}^{'}
\times_{{\mathcal Y}}{\mathcal Y}^{'}$ we take the smallest closed
analytic subspace ${\mathcal X}^{**}$ containing ${\mathcal X}^{*} -
\cup \Pi^{*-1}(p_{j})$ where $\Pi^{*}$ is the composite map
${\mathcal X}^{*} \rightarrow {\mathcal Y}^{'} \rightarrow {\mathcal
Y}$.  Then ${\mathcal X}^{**} \rightarrow {\mathcal Y}^{'}$ is flat.
Finally, resolving the singularities of ${\mathcal X}^{**}$ and of
the components of the singular fiber, we get the required
degeneration ${\mathcal X}^{''}/\Delta^{'}$. \hfill $\Box$\\ 

\begin{corollary} If ${\mathcal X}/\Delta$ is weakly projective, then
each component of ${\mathcal X}_{0}$ is either a fibration over a curve
of genus $\le g(S)$ with fiber genus $\le g(F)$ or a ruled or
rational surface.
\end{corollary}

\noindent {\it Proof :}
Since ${\mathcal X}$ is weakly projective, so is ${\mathcal X}^{''}$ of
Theorem 2 and we apply ([5], Corollary 3.1.4). \hfill $\Box$\\

Next we address to the question of when the hypothesis of Theorem 1
are satisfied for a degeneration ${\mathcal X} \rightarrow \Delta$. As
to the first condition we have\\ a) if $X \rightarrow S$ is a smooth
fibration with both fiber genus and base genus $\ge 2$, then any
deformation of $X$ admits a fibration of the same type ([3], Lemma
7.1), and\\ b) if the fibration is a consequence of a relation among
some deformation invariants, then trivially the first condition of
the hypothesis holds. One notable example of this case is
degenerations of minimal surfaces with $K^{2} < 3 \chi, q \ge 2$.
With this inequality satisfied, the given surface admits a fibration
of fiber genus 2 or 3, the base curve being the image of the
albanese map ([1], Theorem 2.6).\\ The condition on monodromy is
trivially satisfied if $h^{1}(\vert \Gamma \vert) = 0$ where
$\Gamma$ is the dual graph of the singular fiber. For more general
degenerations we have\\

\begin{lemma} In the following cases the hypothesis on monodromy
is satisfied :\\
(a) $X \rightarrow S$ is a smooth fibration with $g(S), g(F) \ge 2$,\\ 
(b) $X$ is minimal and $K^{2} < 3\chi(X), q(X) \ge 2$.
\end{lemma} 

\noindent {\it Proof :}
(a) By ([3], Lemma 7.1), we have a deformation $\Phi : {\mathcal S}
\rightarrow \Delta^{*}$ of $S$, varying continuously with $t \in
\Delta^{*}$, such that $\Pi\vert_{\Delta^{*}}$ factors through
$\Phi$. Therefore, we have an exact sequence $0 \rightarrow
R^{1}\Phi_{*}(\mathbb{C}) \rightarrow R^{1}\Pi_{*}(\mathbb{C})$.  Hence,
the correspondence between representations of $\pi_{1}(\Delta^{*})$
and flat vector bundles on $\Delta^{*}$ shows that, for each $t \in
\Delta^{*}$, $H^{1}(S_{t}, \mathbb{C})$ is an invariant subspace of
$H^{1}({\mathcal X}_{t}, \mathbb{C})$ under the monodromy action.\\ 

(b) In this case the fibration $\Psi_{t}$ being the albanese
fibration we have $H^{1}({\mathcal X}_{t}, \mathbb{C}) =
H^{1}(S_{t}, \mathbb{C})$ and the conclusion follows trivially.
\hfill $\Box$\\


\end{document}